\documentclass[11pt, a4paper, reqno]{amsart}
\usepackage{amssymb, array, amsmath, amscd, pdfpages, enumerate, amsthm, fixltx2e, setspace, hyperref,mathtools}
\pdfoutput=1
\usepackage[latin1]{inputenc} %oli:[utf8]
\usepackage{amssymb}
\usepackage{bm}
\usepackage{graphicx}

\newcommand*{\R}{{\mathbb{R}}}

\renewcommand{\ker}{{\mathcal{N}}}

\newcommand*{\Abs}[2][default]{\ifthenelse{\equal{#1}{default}}{\left\lvert#2\right\rvert}{\ldelim{#1}{\lvert}#2\rdelim{#1}{\rvert}}}
\newcommand*{\Norm}[2][default]{\ifthenelse{\equal{#1}{default}}{\left\lVert#2\right\rVert}{\ldelim{#1}{\lVert}#2\rdelim{#1}{\rVert}}}

\newcommand*{\Iprod}[3][default]{\ifthenelse{\equal{#1}{default}}{\left\langle#2,#3\right\rangle}{\ldelim{#1}{\langle}#2,#3\rdelim{#1}{\rangle}}}
\newcommand*{\Dualpair}[3][default]{\ifthenelse{\equal{#1}{default}}{\left\langle#2,#3\right\rangle}{\ldelim{#1}{\langle}#2,#3\rdelim{#1}{\rangle}}}

\newcommand*{\ddb}[2][1]{\ifthenelse{\equal{#1}{1}}{\frac{d}{d#2}}{\frac{d^{#1}}{d#2^{#1}}}}
\newcommand*{\pd}[3][1]{\ifthenelse{\equal{#1}{1}}{\frac{\partial{#2}}{\partial{#3}}}{\frac{\partial^{#1}{#2}}{\partial#3^{#1}}}}

\newcommand*\lenv{{\hbox{\raisebox{-.15ex}{\rotatebox[origin=c]{50}{$\smallsmile$}}\kern-8.65pt\rotatebox[origin=c]{-25}{$\smallsetminus$}}}}
\newcommand*\uenv{{\hbox{\raisebox{-.0ex}{\rotatebox[origin=c]{-45}{$\smallfrown$}}\kern-5.6pt\raisebox{.2ex}{\rotatebox[origin=c]{-5}{\scriptsize\slash}}}}\,\kern+1.5pt}

\newcommand*{\lupp}[1]{\kern0pt^{\kern0pt u \kern0pt}\kern1pt#1}
\newcommand*{\llow}[1]{\kern0pt^{\kern0pt l \kern0pt}\kern1pt#1}
\newcommand*{\rupp}[1]{#1\kern1pt^{\kern0pt u \kern0pt}\kern0pt}
\newcommand*{\rlow}[1]{#1\kern1pt^{\kern0pt l \kern0pt}\kern0pt}
\newcommand*{\ul}[1]{\kern0pt^{\kern0pt u \kern0pt}\kern0pt#1\kern1pt^{\kern0pt l \kern0pt}\kern0pt}
\newcommand*{\lu}[1]{\kern0pt^{\kern0pt l \kern0pt}\kern0pt#1\kern0pt^{\kern0pt u \kern0pt}\kern0pt}

\DeclareMathOperator{\mmin}{Min}
\DeclareMathOperator{\mmax}{Max}
\newcommand{\mini}[1]{\mmin [\kern1pt #1 \kern1pt ]}
\newcommand{\maxi}[1]{\mmax [\kern1pt #1 \kern1pt ]}

\newcommand{\sa}{\mathrm{sa}}
\newcommand{\Cstar}{$C^*$}

\newcommand*{\bproofname}{Proof}
\newenvironment{bproof}[1][\bproofname]{\begin{proof}[#1]}{\end{proof}}

\newtheorem{thm}{Theorem}[section]
\newtheorem{prop}[thm]{Proposition}
\newtheorem{lemma}[thm]{Lemma}
\newtheorem{cor}[thm]{Corollary}
\theoremstyle{definition}

\newtheorem{remark}[thm]{Remark}
\newtheorem{example}[thm]{Example}

\numberwithin{equation}{section}

\begin{document}

\title[Disjointness, Orthogonality and Commutativity in $C^*$-Algebras]{Disjointness, Orthogonality and Commutativity in $C^*$-Algebras}

\thispagestyle{plain}

\author{Jani Jokela}
\address[J. Jokela]{Mathematics Research Centre, Tampere University, P.O. Box 692, 33101 Tampere, Finland.}
\email{jani.jokela@tuni.fi}

\begin{abstract}
We characterize order-theoretic disjointness of two elements $a$ and $b$ in the self-adjoint part of a $C^*$-algebra $A$ via the condition $aAb = \{0\}$ and the trivial intersection of closed two-sided ideals generated by $a$ and $b$. As a consequence, we show that disjoint complements coincide with the self-adjoint parts of annihilator ideals. We also prove that $A$ is commutative whenever every pair of algebraically orthogonal positive contractions in $A$ is disjoint. Finally, we establish corresponding commutativity criteria involving projections for real rank zero $C^*$-algebras. 
\end{abstract}

\subjclass[2020]{Primary 46L05; Secondary 46B40, 46A40, 06F20}
\keywords{$C^*$-algebra, disjointness, pre-Riesz space, commutativity, orthogonality, annihilator ideal, real rank zero}

\maketitle

\section{Introduction}\label{sec:intro}

We begin by presenting the most important notations and basic facts on $C^*$-algebras used in this paper. 
Throughout this paper, $A$ denotes a $C^*$-algebra, not necessarily unital. We write \[ A_{\sa}=\{x\in A:x^*=x\} \] for its self-adjoint part and \[ A_+=\{c^*c:c\in A\} \] for its positive cone \cite[Theorem~1.3.3]{Pedersen2018}. The closed, convex, and proper real cone \(A_+\) induces a partial order on \(A_{\sa}\), defined by \[ x\leq y \quad\Longleftrightarrow\quad y-x\in A_+. \] With this order, \(A_{\sa}\) is an Archimedean directed partially ordered real Banach space. For \(a,b\in A_{\sa}\) with \(a\leq b\), we denote the corresponding order interval by \[ [a,b]=\{x\in A_{\sa}:a\leq x\leq b\}. \] 
We write \[ \operatorname{Proj}(A)=\{p\in A:p=p^*=p^2\} \] for the set of projections in \(A\), and 
\[ Z(A)=\{z\in A:za=az\text{ for every }a\in A\} \] for the center of \(A\). When $A$ is unital, we denote its identity by $\mathbf{1}$; when $A$ is non-unital, $\widetilde{A}$ denotes its minimal unitization with identity $\mathbf{1}$. %\cite[Proposition~1.1.3]{Pedersen2018}. 
For a Hilbert space \(H\), we denote by \(\mathbf{B}(H)\) the \(C^*\)-algebra of all bounded linear operators on \(H\).

If $A$ is non-commutative then $A_{\sa}$ is one of the central examples of a non-lattice ordered Banach space.  
To extend lattice-theoretic concepts to more general partially ordered vector spaces, van Haandel \cite{vanHaandel1993} introduced the theory of \emph{pre-Riesz spaces}. Attempts to extend the classical notion of disjointness from vector lattices (Riesz spaces) to more general ordered vector spaces have led to several distinct concepts. We will discuss two of these generalizations. For a subset $M$ of a partially ordered vector space $X$, we denote by $M^{u} = \{w \in X : w \ge m \text{ for all } m \in M\}$ the set of all upper bounds of $M$ in $X$.  Two elements $u, v \in X$ are called \emph{disjoint}, written $u \perp v$, if $\{u + v, -(u + v)\}^u =\{u - v, v - u\}^u$. 
Two positive elements $u, v \in X_+$ are called \emph{$D$-disjoint} if $[0, u] \cap [0, v] = \{0\}$. In any vector lattice, both of these conditions coincide with $|u| \wedge |v| = 0$. In any pre-Riesz space, disjointness implies $D$-disjointness for positive elements. The converse fails in general, and $X$ is said to be \emph{weakly pervasive} if every pair of $D$-disjoint elements in $X_+$ is disjoint. %\cite{KV}.

Since $A_{\sa}$ is directed and Archimedean (and hence a pre-Riesz space) for every \Cstar-algebra $A$, it is a natural and fundamental problem to determine what disjointness $u \perp v$ means in terms of the algebraic and ideal-theoretic structure of $A$. 
The purpose of this paper is to give a complete algebraic and ideal-theoretic characterization of disjointness in arbitrary \Cstar-algebras and to explore its structural consequences. Our main result (Theorem~\ref{thm:main}) establishes that for $a, b \in A_{\sa}$, the order-theoretic condition $a \perp b$ is equivalent to the strong algebraic orthogonality condition $aAb = \{0\}$, as well as to the trivial intersection $\overline{AaA} \cap \overline{AbA} = \{0\}$ of the closed two-sided ideals generated by $a$ and $b$. 
This characterization reveals a distinction between disjointness and weaker notions of orthogonality in noncommutative \Cstar-algebras. While $a \perp b$ always implies algebraic orthogonality $ab = 0$ (Corollary~\ref{cor:ortho}), and for $a,b \in A_+$ this in turn implies $D$-disjointness $[0, a] \cap [0, b] = \{0\}$ (Proposition~\ref{prop:ab_zero_interval}), neither implication is reversible in general, as shown by Example~\ref{ex:m2_hierarchy}. Furthermore, if $M\subseteq A_{\sa}$, then applying our main characterization to the disjoint complement $M^d = \{x \in A_{\sa} : x \perp m \text{ for all } m \in M\}$ shows that $M^d = I_M^\perp \cap A_{\sa}$, where $I_M^\perp$ is the annihilator ideal of the closed two-sided ideal generated by $M$ (Corollary~\ref{cor:bands}). 

Commutative $C^*$-algebras have been characterized in terms of their order-theoretic properties. A classical theorem of Sherman \cite{Sherman51} states that $A_{\sa}$ is a vector lattice if and only if $A$ is commutative.  A recent result shows that a $C^*$-algebra $A$ whose self-adjoint part $A_{\sa}$ is a weakly pervasive pre-Riesz space must be commutative \cite[Corollary 5.3]{Jokela26} (see also \cite[Theorem 3.39]{vangaans26} for an analogous result in $JB$-algebras). In this paper, we show that the full strength of weak pervasiveness is not required by proving that commutativity is already forced by apparently weaker disjointness hypotheses. Specifically, we show that $A$ is commutative whenever $ab = 0 \implies a \perp b$ holds merely for pairs of positive contractions in $A$ (Corollary~\ref{cor:comm}). When restricted to projections, the equivalence of $pq = 0$ and $p \perp q$ characterizes centrality of all projections in unital or real rank zero \Cstar-algebras, and hence characterizes commutativity in the real rank zero setting (Corollary~\ref{cor:proj}).

\section{Disjointness in $C^*$-algebras}

We begin with elementary properties of disjoint positive elements that will be used repeatedly.

\begin{lemma}\label{lem:basic} 
Let $A$ be a \Cstar-algebra. 
\begin{enumerate} 
\item[\textup{(a)}] If $a,b\in A_+$ satisfy $a\perp b$, then $x\perp y$ for all $x,y\in A_+$ such that $0\leq x\leq a$ and $0\leq y\leq b$. 
\item[\textup{(b)}] If $u,v\in A_{\sa}$ satisfy $u\perp v$, then $su\perp tv$ for all $s,t\in\R$. 
\end{enumerate} 
\end{lemma}

\begin{bproof} 
\textup{(a)} We use the following immediate consequence of the definition: if $c,d\in A_+$, then $c\perp d$ if and only if every $w\in A_{\sa}$ satisfying $w\geq\pm(c-d)$ also satisfies $w\geq c+d$. Suppose that $a\perp b$, and let $x\in A_+$ satisfy $0\leq x\leq a$. We first show that $x\perp b$. Let $w\in A_{\sa}$ satisfy $w\geq\pm(x-b)$. Then 
\[ w+(a-x) \geq (x-b)+(a-x) = a-b. \] 
Moreover, since $a-x\geq0$, we have 
\[ w+(a-x) \geq w \geq b-x \geq b-a. \] 
Thus $w+(a-x)\geq\pm(a-b)$. Since $a\perp b$, the observation above gives \[ w+(a-x)\geq a+b. \] Subtracting $a-x$ yields $w\geq x+b$. Hence $x\perp b$. Now let $y\in A_+$ satisfy $0\leq y\leq b$. %By symmetry, $b\perp x$. 
Applying the preceding argument to $x\perp b$ gives $x\perp y$. 

\par\smallskip \noindent 
\textup{(b)} For a subset $S\subseteq A_{\sa}$, let \[ S^{\mathrm d} = \{x\in A_{\sa}:x\perp y \text{ for every }y\in S\} \] denote its disjoint complement. Since $A_{\sa}$ is a pre-Riesz space, the disjoint complement of every subset is a linear subspace by \cite[Proposition 4.1.5]{KV}. Suppose that $u\perp v$. Then $v\in\{u\}^{\mathrm d}$, and hence $tv\in\{u\}^{\mathrm d}$ for every $t\in\mathbb{R}$. Thus $u\perp tv$. But then $u\in\{tv\}^{\mathrm d}$, and since $\{tv\}^{\mathrm d}$ is a linear subspace, $su\in\{tv\}^{\mathrm d}$ for every $s\in\mathbb{R}$. Therefore, $su\perp tv$. 
\end{bproof}

To prove our main results, we need a couple of technical lemmas on $C^*$-algebras. 
%Throughout this paper, all representations of $C^*$-algebras are assumed to be nondegenerate.

\begin{lemma}\label{lem:irred_compression}
Let $A$ be a \Cstar-algebra, let $\pi \colon A \to \mathbf{B}(H)$ be an irreducible representation, and let $K \subseteq H$ be a 
closed subspace with orthogonal projection $P_K \in \mathbf{B}(H)$. If $a, b \in A_{\sa}$ satisfy $a \perp b$ in $A_{\sa}$, then
\begin{enumerate}
    \item[\textup{(a)}] $\pi(a) \perp \pi(b)$ in $\mathbf{B}(H)_{\sa}$,
    \item[\textup{(b)}] $(P_K \pi(a) P_K)|_K \perp (P_K \pi(b) P_K)|_K$ in $\mathbf{B}(K)_{\sa}$.
\end{enumerate}
\end{lemma}

\begin{bproof}
\textup{(a)} Set $\widetilde{a}=\pi(a)$ and $\widetilde{b}=\pi(b)$. We prove that 
\[ 
\{\widetilde{a}+\widetilde{b},-(\widetilde{a}+\widetilde{b})\}^u = \{\widetilde{a}-\widetilde{b},\widetilde{b}-\widetilde{a}\}^u. 
\] 
Let $W\in\mathbf{B}(H)_{\sa}$ satisfy \[ W\geq\pm(\widetilde{a}-\widetilde{b}). \] Since $\pi$ is irreducible, it is nondegenerate and Schur's lemma gives $\pi(A)'=\mathbb{C}I_H$. Hence $\pi(A)''= (\mathbb{C}I_H)'=\mathbf{B}(H)$, and the von Neumann bicommutant theorem \cite[Corollary~2.2.5]{Pedersen2018} implies that the $^*$-subalgebra $\pi(A)$ is strongly dense in $\mathbf{B}(H)$. By Kaplansky's density theorem \cite[Theorem~2.3.3]{Pedersen2018}, there exists a net $(T_\alpha)$ in $\pi(A)_{\sa}$ such that 
\[ 
\|T_\alpha\|\leq\|W\| \quad\text{for every $\alpha$ and} \quad T_\alpha\longrightarrow W \quad\text{strongly}. 
\] 
For each $\alpha$, choose $z_\alpha\in A$ such that $\pi(z_\alpha)=T_\alpha$. There is no guarantee that each $z_{\alpha}$ is in $A_{\sa}$, so we define 
\[ 
w_\alpha=\frac{z_\alpha+z_\alpha^*}{2}\in A_{\sa}. 
\] 
Since $T_\alpha$ is self-adjoint, we have 
\[ 
\pi(w_\alpha) = \frac{T_\alpha+T_\alpha^*}{2} = T_\alpha. 
\] 
Consequently, 
\[ 
\|\pi(w_\alpha)\|\leq\|W\| \quad\text{and}\quad \pi(w_\alpha)\longrightarrow W \quad\text{strongly}. 
\]  
For each $\alpha$, define the elements 
\[ 
e_\alpha = \bigl((a-b)-w_\alpha\bigr)_+ + \bigl((b-a)-w_\alpha\bigr)_+ \in A_+ 
\] 
and set $\widehat{w}_\alpha=w_\alpha+e_\alpha\in A_{\sa}$.  
Since $u_+,u_-\in A_+$ and $u_+=u+u_-$ for all $u\in A_{\sa}$, we have 
\begin{align*} 
\widehat{w}_\alpha &\geq w_\alpha+\bigl((a-b)-w_\alpha\bigr)_+ 
\\ &= w_\alpha+\bigl((a-b)-w_\alpha\bigr) +\bigl((a-b)-w_\alpha\bigr)_- 
\\ &= a-b+\bigl((a-b)-w_\alpha\bigr)_- 
\\ &\geq a-b. 
\end{align*} 
Similarly, $\widehat{w}_\alpha\geq b-a$. Thus $\widehat{w}_\alpha \in \{a-b,b-a\}^u$, and since $a\perp b$, we have
\[ \{a-b,b-a\}^u = \{a+b,-(a+b)\}^u. \] It follows that $\widehat{w}_\alpha\geq\pm(a+b)$, and since representations are order-preserving, we have $\pi(\widehat{w}_\alpha) \geq \pm(\widetilde{a}+\widetilde{b})$.  Because $(\pi(w_\alpha))$ is norm-bounded and converges strongly to $W$, the nets 
\[ 
(\widetilde{a}-\widetilde{b})-\pi(w_\alpha) \quad\text{and}\quad (\widetilde{b}-\widetilde{a})-\pi(w_\alpha) 
\] 
are norm-bounded and converge strongly to \[ (\widetilde{a}-\widetilde{b})-W \quad\text{and}\quad (\widetilde{b}-\widetilde{a})-W, \] respectively. The assumption $W\geq\pm(\widetilde{a}-\widetilde{b})$ implies that both limit operators are negative, and so their positive parts are zero. Therefore 
%%Since $*$-homomorphisms preserve continuous functional calculus, we have 
\begin{align*} 
\pi(e_\alpha) &= \pi\left(\bigl((a-b)-w_\alpha\bigr)_+\right) + \pi\left(\bigl((b-a)-w_\alpha\bigr)_+\right) 
\\ &= \bigl((\widetilde{a}-\widetilde{b})-\pi(w_\alpha)\bigr)_+ 
 + \bigl((\widetilde{b}-\widetilde{a})-\pi(w_\alpha)\bigr)_+ \longrightarrow 0 
\end{align*} 
strongly. Consequently, 
%%By the strong continuity of continuous functional calculus on norm-bounded subsets of $\mathbf{B}(H)_{\sa}$, 
%%\cite[Proposition~2.3.2]{Pedersen2018}, tuo viite on väärin
%%it follows that \[ \pi(e_\alpha)\longrightarrow 0 \quad\text{strongly}. \] Therefore, 
\[ \pi(\widehat{w}_\alpha) = \pi(w_\alpha)+\pi(e_\alpha) \longrightarrow W \quad\text{strongly}. \] 
The positive cone $\mathbf{B}(H)_+$ is strongly closed, so passing to the strong limit in 
$\pi(\widehat{w}_\alpha) \geq \pm(\widetilde{a}+\widetilde{b})$  
gives $W\geq\pm(\widetilde{a}+\widetilde{b})$.  
We have proved that 
\[ 
\{\widetilde{a}-\widetilde{b}, \widetilde{b}-\widetilde{a}\}^u \subseteq \{\widetilde{a}+\widetilde{b}, -(\widetilde{a}+\widetilde{b})\}^u. 
\] 
To prove the reverse inclusion, observe that $a\perp-b$. Indeed, replacing $b$ by $-b$ merely interchanges the two upper-bound sets in the definition of disjointness. Applying the inclusion proved above to the disjoint pair $a,-b$ gives 
\[ 
\{\widetilde{a}+\widetilde{b}, -(\widetilde{a}+\widetilde{b})\}^u \subseteq \{\widetilde{a}-\widetilde{b}, \widetilde{b}-\widetilde{a}\}^u. 
\] 
The two sets of upper-bounds are therefore equal, and hence $\widetilde{a}\perp\widetilde{b}$ in $\mathbf{B}(H)_{\sa}$.

\par\smallskip \noindent
\textup{(b)} Set $A_K=(P_K\widetilde{a}P_K)|_K$ and $B_K=(P_K\widetilde{b}P_K)|_K$. Then $A_K,B_K\in\mathbf{B}(K)_{\sa}$. We first prove that 
\[
\{A_K-B_K,B_K-A_K\}^u \subseteq \{A_K+B_K,-(A_K+B_K)\}^u.
\]
Let $W_1\in\mathbf{B}(K)_{\sa}$ satisfy $W_1\geq\pm(A_K-B_K)$. With respect to the orthogonal decomposition $H=K\oplus K^\perp$, write 
\[ 
\widetilde{a}-\widetilde{b} = \begin{pmatrix} A_K-B_K & R\\ R^* & C \end{pmatrix}, 
\] 
where $R=P_K(\widetilde{a}-\widetilde{b})|_{K^\perp} \in\mathbf{B}(K^\perp,K)$ and $C=((\mathbf{1}-P_K)(\widetilde{a}-\widetilde{b}) (\mathbf{1}-P_K))|_{K^\perp} \in\mathbf{B}(K^\perp)_{\sa}$. 

Fix $\varepsilon>0$, choose $M>\|C\|+\varepsilon^{-1}\|R\|^2$, and define $W_{\varepsilon,M} =(W_1+\varepsilon\mathbf{1}_K)\oplus M\mathbf{1}_{K^\perp}\in\mathbf{B}(H)_{\sa}$. We will show that $W_{\varepsilon,M}\mp(\widetilde{a}-\widetilde{b})$ is positive semidefinite. Let $z\in H$. Then $z=x+y$, where $x\in K$ and $y\in K^\perp$, and we have 
\begin{align*} 
&\left\langle \bigl(W_{\varepsilon,M} \mp(\widetilde{a}-\widetilde{b})\bigr)(x+y),x+y \right\rangle \\ 
&\quad= \left\langle \bigl(W_1\mp(A_K-B_K)+\varepsilon\mathbf{1}_K\bigr)x,x \right\rangle \mp2\operatorname{Re}\langle Ry,x\rangle +\left\langle \bigl(M\mathbf{1}_{K^\perp}\mp C\bigr)y,y \right\rangle \\ 
&\quad\geq \varepsilon\|x\|^2-2\|R\|\|x\|\|y\| +(M-\|C\|)\|y\|^2 \\ 
&\quad= \bigl(\varepsilon^{1/2}\|x\| -\varepsilon^{-1/2}\|R\|\|y\|\bigr)^2 +\bigl(M-\|C\|-\varepsilon^{-1}\|R\|^2\bigr)\|y\|^2 \geq0. 
\end{align*} 
Here we used $W_1\mp(A_K-B_K)\geq 0$ and Cauchy--Schwarz to get $|\operatorname{Re}\langle Ry,x\rangle| \leq\|R\|\|x\|\|y\|$. It follows that $W_{\varepsilon,M}\geq\pm(\widetilde{a}-\widetilde{b})$. By part~\textup{(i)}, we have $\widetilde{a}\perp\widetilde{b}$ in $\mathbf{B}(H)_{\sa}$. Hence $W_{\varepsilon,M}\geq\pm(\widetilde{a}+\widetilde{b})$. Compressing these inequalities to $K$ gives $W_1+\varepsilon\mathbf{1}_K\geq\pm(A_K+B_K)$. Since $\mathbf{B}(K)_+$ is norm closed, letting $\varepsilon\to0^+$ yields $W_1\geq\pm(A_K+B_K)$. We have therefore proved that 
\[ 
\{A_K-B_K,B_K-A_K\}^u \subseteq \{A_K+B_K,-(A_K+B_K)\}^u. 
\] 
To prove the reverse inclusion, observe that $a\perp-b$, since replacing $b$ by $-b$ interchanges the two upper-bound sets in the definition of disjointness. By part~\textup{(a)}, it follows that $\widetilde{a}\perp-\widetilde{b}$. Applying the inclusion proved above to the pair $\widetilde{a},-\widetilde{b}$ gives 
\[ 
\{A_K+B_K,-(A_K+B_K)\}^u \subseteq \{A_K-B_K,B_K-A_K\}^u. 
\] 
Thus the two upper-bound sets are equal, and consequently $A_K\perp B_K$ in $\mathbf{B}(K)_{\sa}$. 
\end{bproof}

%_____________________________

\begin{lemma}\label{lem:low_dim_matrices}
Let $d \in \{1, 2\}$. If $U, V \in M_d(\mathbb{C})_{\sa}$ satisfy $U \perp V$ in $M_d(\mathbb{C})_{\sa}$, then $U = 0$ or $V = 0$.
\end{lemma}

\begin{bproof} 
Suppose that $U\neq 0$ and $V\neq 0$. Since $U\perp V$, we have \[ \{U+V,-(U+V)\}^u = \{U-V,V-U\}^u. \] Now $|U-V|\geq\pm(U-V)$, so $|U-V|$ is an upper bound of $\{U-V,V-U\}$. Hence it is also an upper bound of $\{U+V,-(U+V)\}$, and therefore $|U-V|\geq\pm(U+V)$. Similarly, since $|U+V|\geq\pm(U+V)$, we obtain $|U+V|\geq\pm(U-V)$. We first record a consequence of these inequalities. Suppose that $X,Y\in M_d(\mathbb{C})_{\sa}$ satisfy $X\geq\pm Y$. Then $X-Y\geq0$ and $X+Y\geq0$, so \[ (X-Y)^{1/2}(X+Y)(X-Y)^{1/2}\geq0. \] Using the linearity of the trace and the identity $\operatorname{Tr}(XY)=\operatorname{Tr}(YX)$, we obtain 
\begin{align*} 
0 &\leq \operatorname{Tr}\bigl( (X-Y)^{1/2}(X+Y)(X-Y)^{1/2} \bigr) \\ 
&= \operatorname{Tr}\bigl((X-Y)(X+Y)\bigr) \\ 
&= \operatorname{Tr}(X^2)-\operatorname{Tr}(Y^2). 
\end{align*} 
Applying this first with $X=|U-V|$ and $Y=U+V$ gives 
\[ 
\operatorname{Tr}\bigl((U-V)^2\bigr) \geq \operatorname{Tr}\bigl((U+V)^2\bigr). 
\] 
Applying it with $X=|U+V|$ and $Y=U-V$ gives the reverse inequality. Consequently, 
\[ 
\operatorname{Tr}\bigl((U-V)^2\bigr) = \operatorname{Tr}\bigl((U+V)^2\bigr). 
\] 
Set $X=|U-V|$, $Y=U+V$, $P=X-Y$, and $Q=X+Y$. The inequality $|U-V|\geq\pm(U+V)$ shows that $P,Q\geq0$. Moreover, 
\begin{align*} 
\operatorname{Tr}(P^{1/2}QP^{1/2}) &= \operatorname{Tr}(PQ) \\ 
&= \operatorname{Tr}\bigl((X-Y)(X+Y)\bigr) \\ 
&= \operatorname{Tr}(X^2)-\operatorname{Tr}(Y^2) \\ 
&= \operatorname{Tr}\bigl((U-V)^2\bigr) - \operatorname{Tr}\bigl((U+V)^2\bigr) =0. 
\end{align*} 
Since $P^{1/2}QP^{1/2}\geq0$ and a positive semidefinite matrix with trace zero must be zero, it follows that $P^{1/2}QP^{1/2}=0$. Equivalently, 
\[(Q^{1/2}P^{1/2})^*(Q^{1/2}P^{1/2})=0,  \]
so $Q^{1/2}P^{1/2}=0$. Hence $QP=0$, and therefore $PQ=(QP)^*=0$. It follows that 
\[
0 =\frac{1}{2}(PQ+QP)  =X^2-Y^2 =|U-V|^2-(U+V)^2.
\] 
Since $|U-V|^2=(U-V)^2$, we obtain $(U-V)^2=(U+V)^2$, and hence \[ UV+VU=0. \] If $d=1$, then $U$ and $V$ are real scalars, and the preceding identity becomes $2UV=0$. This contradicts $U\neq0$ and $V\neq0$. Thus $d=2$. Choose an orthonormal eigenbasis of $U$ and write $U=\operatorname{diag}(\lambda_1,\lambda_2)$, where, after possibly interchanging the basis vectors, we may assume that $\lambda_1\neq0$. Writing 
\[
V=\begin{pmatrix}V_{11} & V_{12} \\ V_{21} & V_{22} \end{pmatrix}
\]
in this basis, the identity $UV+VU=0$ becomes 
\[ 
\begin{pmatrix} 2\lambda_1V_{11} & (\lambda_1+\lambda_2)V_{12} \\ (\lambda_1+\lambda_2)V_{21} & 2\lambda_2V_{22} \end{pmatrix} =0. 
\] 
Thus $V_{11}=0$, $(\lambda_1+\lambda_2)V_{12}=0$, and $\lambda_2V_{22}=0$. Since $V$ is self-adjoint, $V_{21}=\overline{V_{12}}$. As $V\neq0$, we have two possible cases: either $\lambda_2=0$, in which case we have $V_{12}=0$ and $V_{22}\neq 0$; or $\lambda_2\neq 0$, in which case we must have $V_{22}=0$ and $V_{12}\neq 0$, and $(\lambda_1+\lambda_2)V_{12}=0$ implies that $\lambda_2=-\lambda_1$. We consider both cases.

\par\smallskip \noindent
\textbf{Case 1:} If $V_{12}=0$ and $V_{22}\neq 0$,   
then $U=\operatorname{diag}(\lambda_1,0)$ and $V=\operatorname{diag}(0,V_{22})$. Disjointness is preserved under independent sign changes. Thus, replacing $U$ by $\operatorname{sgn}(\lambda_1)U$ and $V$ by $\operatorname{sgn}(V_{22})V$, we may assume that 
\[ 
U= \begin{pmatrix} \alpha&0\\ 0&0 \end{pmatrix}, \qquad V= \begin{pmatrix} 0&0\\ 0&\beta \end{pmatrix}, 
\] 
where $\alpha,\beta>0$. Define 
\[
W= \begin{pmatrix} 2\alpha&\sqrt{3\alpha\beta}\\ \sqrt{3\alpha\beta}&2\beta \end{pmatrix}. 
\] 
Then 
\[ 
W-(U-V) = \begin{pmatrix} \alpha&\sqrt{3\alpha\beta}\\ \sqrt{3\alpha\beta}&3\beta \end{pmatrix} \quad\text{and}\quad W-(V-U) = \begin{pmatrix} 3\alpha&\sqrt{3\alpha\beta}\\ \sqrt{3\alpha\beta}&\beta \end{pmatrix}. 
\] 
Both matrices have positive diagonal entries and determinant zero, and are therefore positive semidefinite. Thus $W\geq\pm(U-V)$. On the other hand, 
\[ 
W-(U+V) = \begin{pmatrix} \alpha&\sqrt{3\alpha\beta}\\ \sqrt{3\alpha\beta}&\beta \end{pmatrix} 
\] 
has determinant $-2\alpha\beta<0$, and hence is not positive semidefinite. Therefore $W\not\geq U+V$. 
This contradicts the assumption that $U\perp V$.

\par\smallskip \noindent
\textbf{Case 2:} If $V_{22}=0$, $V_{12}\neq 0$ and $\lambda_2=-\lambda_1$,   
then we write $V_{12}=re^{i\theta}$, where $r>0$ and $\theta\in\mathbb{R}$, and set $D=\operatorname{diag}(1,e^{-i\theta})$. Then $U=D^*UD$ and 
\[ 
D^*VD= \begin{pmatrix} 0&r\\ r&0 \end{pmatrix}. 
\] 
Since $D$ is unitary, the mappings $X\mapsto D^*XD$ and $X\mapsto cX$ (for $c>0$) are order isomorphisms and hence they preserve disjointness. Therefore, replacing $U$ by $\operatorname{sgn}(\lambda_1)U$ and multiplying both $U$ and $D^*VD$ by the common positive scalar $|\lambda_1|^{-1}$, we may assume that 
\[ 
U= \begin{pmatrix} 1&0\\ 0&-1 \end{pmatrix}, \qquad V= \begin{pmatrix} 0&t\\ t&0 \end{pmatrix}, 
\] 
where $t=r/|\lambda_1|>0$. Choose $\delta>0$ such that $\delta(2+\delta)<4t^2$, and define 
\[ 
W_1= \begin{pmatrix} 1+\delta&-t\\ -t&\dfrac{4t^2}{2+\delta}+1 \end{pmatrix}. 
\] 
We have 
\[ 
W_1-(U-V) = \begin{pmatrix} \delta&0\\ 0&\dfrac{4t^2}{2+\delta}+2 \end{pmatrix} \geq0. 
\] 
Moreover, 
\[ 
W_1+(U-V) = \begin{pmatrix} 2+\delta&-2t\\ -2t&\dfrac{4t^2}{2+\delta} \end{pmatrix}. 
\] 
This matrix has positive diagonal entries and determinant zero, and is therefore positive semidefinite. Consequently, $W_1\geq\pm(U-V)$. However, 
\[ 
W_1-(U+V) = \begin{pmatrix} \delta&-2t\\ -2t&\dfrac{4t^2}{2+\delta}+2 \end{pmatrix} 
\] 
has determinant \[ \frac{2\delta(2+\delta)-8t^2}{2+\delta}<0. \] Thus $W_1\not\geq U+V$, and this again contradicts the assumption $U\perp V$.  
Both cases lead to a contradiction, and hence $U=0$ or $V=0$. 
\end{bproof}

%___________________________MAIN THEOREM__________________________

We are now ready to state and prove our main theorem.

\begin{thm}\label{thm:main}%[Characterization of Disjointness in \Cstar-Algebras]
Let $A$ be a \Cstar-algebra and let $u, v \in A_{\sa}$. The following conditions are equivalent:
\begin{enumerate}[(a)]
    \item $u \perp v$ in $A_{\sa}$.
    \item $uAv = \{0\}$ (that is, $uxv = 0$ for every $x \in A$).
    \item The closed two-sided ideals generated by $u$ and $v$ have trivial intersection: $\overline{AuA} \cap \overline{AvA} = \{0\}$.
\end{enumerate}
\end{thm}

\begin{bproof} 
$\textup{(b)}\iff\textup{(c)}$: Let $AuA$ and $AvA$ denote the linear spans of $\{xuy : x,y\in A\}$ and $\{xvy : x,y\in A\}$, respectively. Since $A$ admits an approximate unit, the closures $I_u=\overline{AuA}$ and $I_v=\overline{AvA}$ are the closed two-sided ideals generated by $u$ and $v$, respectively. 
Suppose that $uAv=\{0\}$. Then products of elements of $AuA$ and $AvA$ vanish, and hence $I_uI_v=\{0\}$. If $z\in I_u\cap I_v$, then $z^*\in I_u$ and $z\in I_v$, so $z^*z=0$ and therefore $z=0$. 

Conversely, if $I_u\cap I_v=\{0\}$, then $uxv\in I_u\cap I_v$ for every $x\in A$, and hence $uxv=0$.

\par\smallskip \noindent 
$\textup{(b)}\implies\textup{(a)}$: Suppose that $uAv=\{0\}$. We identify $A^{**}$ with $\pi_u(A)''\subseteq\mathbf{B}(H_u)$ via the universal representation (see \cite[Section~3.7]{Pedersen2018}). We first show that $uA^{**}v=\{0\}$. Let $y\in A^{**}$. By Kaplansky's density theorem, there exists a norm-bounded net $(x_\alpha)$ in $A$ that converges strongly to $y$. Since $ux_\alpha v=0$ for every $\alpha$ and multiplication on the left and right by fixed operators is strongly continuous, we obtain $uyv = \lim_\alpha ux_\alpha v = 0$ strongly. Thus $uA^{**}v=\{0\}$.

Let $p \in A^{**}$ be the range projection of $u$ (the orthogonal projection onto $\overline{\operatorname{ran}(u)} = \ker(u)^\perp$), so that $pu = up = u$ and $\ker(p) = \ker(u)$.  
Since $u(sv) = 0$ for every $s\in A^{**}$, we have $\operatorname{ran}(sv) \subseteq \ker(u) = \ker(p)$, and therefore $psv = 0$ for all $s\in A^{**}$. Let $z = c(p) \in Z(A^{**})$ be the central cover of the projection $p$; by \cite[2.6.2]{Pedersen2018}, $z$ is a central projection satisfying $p \leq z$, whence $zp = pz = p$ and consequently
\[
zu = zpu = pu = u \quad\text{and}\quad uz = u.
\]
Moreover, by \cite[Lemma 2.6.3]{Pedersen2018}, $z = \bigvee q^*pq$, where $q$ ranges over the unitary group of $A^{**}$. For every unitary $q \in A^{**}$, we have $(q^*pq)v = q^*(pqv) = 0$, so $\operatorname{ran}(v) \subseteq \ker(q^*pq) = \operatorname{ran}(q^*pq)^\perp$. Taking the intersection over all unitary $q \in A^{**}$ yields $\operatorname{ran}(v) \subseteq \ker(z)$, and hence $zv = 0$. Taking adjoints and using $z = z^*$ and $v = v^*$, we also obtain $vz = 0$. It follows that $(\mathbf{1}-z)u=0$ and $(\mathbf{1}-z)v=v$.

Let $w\in A_{\sa}$ satisfy $w\geq\pm(u-v)$. Equivalently, 
\[ w-(u-v)\geq0 \quad\text{and}\quad w+(u-v)\geq0. \] 
Since $z$ is a central projection, the map 
\[ x\longmapsto z^{1/2}xz^{1/2} = zx \] 
is positive on $A^{**}$. Applying this map to $w-(u-v)\geq0$ gives \[ zw-u\geq 0, \] where we used $zu=u$ and $zv=0$. Applying it to $w+(u-v)\geq0$ gives \[ zw+u\geq 0. \] Therefore, $zw\geq\pm u$. Similarly, the map 
\[ x\longmapsto (\mathbf{1}-z)^{1/2} x (\mathbf{1}-z)^{1/2} = (\mathbf{1}-z)x \] 
is positive on $A^{**}$. Applying this map to $w-(u-v)\geq0$ gives 
\[ (\mathbf{1}-z)w+v\geq0, \] 
where we used $(\mathbf{1}-z)u=0$ and $(\mathbf{1}-z)v=v$. Applying it to $w+(u-v)\geq0$ gives 
\[ (\mathbf{1}-z)w-v\geq0. \] 
Thus $(\mathbf{1}-z)w\geq\pm v$. Adding the inequalities with matching signs, we obtain 
\[ w = zw+(\mathbf{1}-z)w \geq u+v \] and \[ w = zw+(\mathbf{1}-z)w \geq -(u+v). \] 
Hence $w\geq\pm(u+v)$ in $A^{**}$. Since $w,u,v\in A_{\sa}$ and \[ A_+ = A\cap(A^{**})_+, \] these inequalities also hold in $A_{\sa}$. We have therefore proved that 
\[ \{u-v,v-u\}^u \subseteq \{u+v,-(u+v)\}^u. \] 
Finally, $uA(-v)=\{0\}$, so applying the preceding argument to the pair $u,-v$ gives 
\[ \{u+v,-(u+v)\}^u \subseteq \{u-v,v-u\}^u. \] 
The two sets of upper-bounds are therefore equal, and so $u\perp v$ in $A_{\sa}$.

\par\smallskip\noindent $\textup{(a)}\implies\textup{(b)}$: Suppose that $u\perp v$ in $A_{\sa}$. Assume to the contrary, that $uAv\neq\{0\}$. Then there exists $x\in A$ such that $uxv\neq0$. Since irreducible representations separate the points of a $C^*$-algebra, 
there exists a nondegenerate irreducible representation $\pi\colon A\to\mathbf{B}(H)$ such that $\pi(uxv) = \pi(u)\pi(x)\pi(v) \neq0$  (see \cite[Theorem 2.7.3]{dix}). 
In particular, $\pi(u)\neq0$ and $\pi(v)\neq0$. Since $\pi(u)$ and $\pi(v)$ are nonzero self-adjoint operators, there exist unit vectors $\xi,\eta\in H$ such that 
\[ \langle\pi(u)\xi,\xi\rangle\neq0 \quad\text{and}\quad \langle\pi(v)\eta,\eta\rangle\neq0. \] 
Indeed, a self-adjoint operator whose quadratic form vanishes on every vector is zero. Let $K=\operatorname{span}\{\xi,\eta\}$, and let $P_K\in\mathbf{B}(H)$ be the orthogonal projection onto $K$. Then $d=\dim K\in\{1,2\}$. Set \[ U_0=(P_K\pi(u)P_K)|_K \quad\text{and}\quad V_0=(P_K\pi(v)P_K)|_K. \] By Lemma~\ref{lem:irred_compression}\textup{(b)}, we have $U_0\perp V_0$ in $\mathbf{B}(K)_{\sa}$. After choosing an orthonormal basis of $K$, we may identify $\mathbf{B}(K)_{\sa}$ with $M_d(\mathbb{C})_{\sa}$. Under this identification, disjointness is preserved. Lemma~\ref{lem:low_dim_matrices} therefore implies that $U_0=0$ or $V_0=0$. On the other hand, since $\xi,\eta\in K$, we have 
\[ \langle U_0\xi,\xi\rangle = \langle\pi(u)\xi,\xi\rangle \neq0 \] and \[ \langle V_0\eta,\eta\rangle = \langle\pi(v)\eta,\eta\rangle \neq0. \] 
Thus both $U_0$ and $V_0$ are nonzero, which is a contradiction. Consequently, $uAv=\{0\}$.
\end{bproof}

%__________korollaarit:______________________

%%For $u\in A_{\sa}$, we define $|u|=(u^2)^{1/2}$.

As an immediate consequence of the ideal-theoretic characterization, disjointness of self-adjoint elements depends only on their $C^*$-algebraic absolute values.

\begin{cor}\label{cor:absolute_values} 
Let $A$ be a \Cstar-algebra, and let $u,v\in A_{\sa}$. Then \[ u\perp v \quad\Longleftrightarrow\quad |u|\perp|v|. \] 
\end{cor} 

\begin{bproof} 
For every $u\in A_{\sa}$, the elements $u$ and $|u|$ generate the same closed two-sided ideal: \[ \overline{AuA} = \overline{A|u|A}. \] Indeed, since $u\in\overline{AuA}$ and $\overline{AuA}$ is a \Cstar-subalgebra of $A$, continuous functional calculus gives $|u|\in\overline{AuA}$. 

Conversely, computations involving the identity are carried out in the minimal unitization $\widetilde A$. For $n\in\mathbb{N}$, set \[ c_n = u(|u|+n^{-1}\mathbf{1})^{-1}. \] Then $c_n\in A$, since it is obtained from $u$ by continuous functional calculus using the function $f_n(t)=t/(|t|+n^{-1})$, which vanishes at zero. Moreover, $c_n|u|\in\overline{A|u|A}$ and \[ \left\| u-c_n|u| \right\| \leq \frac{1}{n}. \] Thus $u\in\overline{A|u|A}$, proving the claimed equality of ideals. Applying Theorem~\ref{thm:main}, we obtain 
\begin{align*} 
u\perp v &\quad\Longleftrightarrow\quad \overline{AuA}\cap\overline{AvA}=\{0\} \\ 
&\quad\Longleftrightarrow\quad \overline{A|u|A}\cap\overline{A|v|A}=\{0\} \\ 
&\quad\Longleftrightarrow\quad |u|\perp|v|. 
\end{align*} 
\end{bproof}

\begin{cor}\label{cor:ortho}%[Algebraic Orthogonality]
Let $A$ be a \Cstar-algebra and let $u,v\in A_{\sa}$. Then $u \perp v$ implies $uv = vu= 0$.
\end{cor}

\begin{bproof}
By Theorem~\ref{thm:main}, $u\perp v$ implies $uAv=\{0\}$. If $(e_\lambda)$ is an approximate unit for $A$, then $ue_\lambda v=0$ for every $\lambda$, and hence $uv=\lim_\lambda ue_\lambda v=0$. Taking adjoints gives $vu=0$.
\end{bproof}

\begin{cor}\label{cor:bands}
Let $A$ be a \Cstar-algebra, and let $M\subseteq A_{\sa}$. Define the disjoint complement of $M$ by 
\[ 
M^d = \{u\in A_{\sa} : u\perp m \text{ for every }m\in M\}, 
\] 
and let 
\[ 
I_M = \overline{\operatorname{span}} \{xmy:x,y\in A,\ m\in M\} 
\] 
be the closed two-sided ideal generated by $M$. Then \[ M^d = I_M^\perp\cap A_{\sa}, \] where 
\[ 
I_M^\perp = \{x\in A:xI_M=I_Mx=\{0\}\} 
\] 
is the annihilator ideal of $I_M$ in $A$. 
\end{cor}

\begin{bproof} 
We first note that $M\subseteq I_M$. Indeed, let $(e_\lambda)$ be an approximate unit for $A$. For every $m\in M$, we have $e_\lambda m e_\lambda\in I_M$ and $e_\lambda m e_\lambda\to m$ in norm. Since $I_M$ is norm closed, it follows that $m\in I_M$. Let $u\in A_{\sa}$. By Theorem~\ref{thm:main}, \[ u\in M^{d} \quad\Longleftrightarrow\quad uAm=\{0\} \text{ for every }m\in M. \] Suppose that $uAm=\{0\}$ for every $m\in M$. For $a,b\in A$ and $m\in M$, we have \[ u(amb) = (uam)b = 0. \] Thus $u$ annihilates $\operatorname{span}\{amb:a,b\in A,\ m\in M\}$. By norm continuity of multiplication, it follows that $uI_M=\{0\}$. 

Conversely, suppose that $uI_M=\{0\}$. Since $M\subseteq I_M$ and $I_M$ is a two-sided ideal, we have $am\in I_M$ for every $a\in A$ and $m\in M$. Consequently, $uam=0$ for every $a\in A$ and $m\in M$. Thus $uAm=\{0\}$ for every $m\in M$. We have therefore shown that \[ u\in M^{d} \quad\Longleftrightarrow\quad uI_M=\{0\}. \] Since $u=u^*$ and $I_M$ is self-adjoint, the condition $uI_M=\{0\}$ is equivalent to $I_Mu=\{0\}$. Indeed, if $uI_M=\{0\}$ and $x\in I_M$, then $x^*\in I_M$, so \[ xu = (ux^*)^* = 0. \] The converse follows by the same argument. Consequently, \[ u\in M^d \quad\Longleftrightarrow\quad uI_M=I_Mu=\{0\} \quad\Longleftrightarrow\quad u\in I_M^\perp\cap A_{\sa}. \] Therefore, \[ M^{d} = I_M^\perp\cap A_{\sa}. \] 
\end{bproof}

%_________________________________________________________

\section{Characterizations of commutativity}

As noted in Section \ref{sec:intro}, the following implications hold for all $a,b\in A_+$:
\[
a \perp b \implies ab= 0 \implies [0,a] \cap [0,b] = \{0\}.
\]
The implication $a\perp b \implies ab=0$ was already proved in Corollary \ref{cor:ortho}. The next result establishes the last implication.

\begin{prop}\label{prop:ab_zero_interval} 
Let $A$ be a \Cstar-algebra, and let $a,b\in A_+$ satisfy $ab=0$. Then $[0,a]\cap[0,b]=\{0\}$. 
\end{prop} 

\begin{bproof} 
Let $c\in[0,a]\cap[0,b]$. Then $c\in A_+$ and $0\leq c\leq a,b$. Conjugating $0\leq c\leq a$ by $b$ gives \[ 0\leq bcb\leq bab. \] Since $ab=0$, we have $bab=b(ab)=0$. Hence $bcb=0$. Let $c^{1/2}\in A_+$ be the positive square root of $c$. Since 
\[ bcb = (c^{1/2}b)^*(c^{1/2}b) = 0, \] 
the \Cstar-identity implies that $c^{1/2}b=0$. Now conjugate the inequality $0\leq c\leq b$ by $c^{1/2}$. This gives 
\[ 0 \leq c^{1/2}cc^{1/2} \leq c^{1/2}bc^{1/2}. \] 
The left-hand side is $c^2$, while the right-hand side vanishes because 
\[ c^{1/2}bc^{1/2} = (c^{1/2}b)c^{1/2} = 0. \] 
Therefore, $0\leq c^2\leq0$, and hence $c^2=0$. Finally, since $c=c^*$, the \Cstar-identity gives $\|c\|^2=\|c^2\|=0$. Thus $c=0$, proving that $[0,a]\cap[0,b]=\{0\}$. 
\end{bproof}

The following example shows that, in general, algebraic orthogonality sits properly between disjointness and $D$-disjointness.

\begin{example}\label{ex:m2_hierarchy}
In a noncommutative \Cstar-algebra, neither implication in
\[
a\perp b \implies ab=0 \implies [0,a]\cap[0,b]=\{0\}
\]
is reversible, even for rank-one projections in $A=M_2(\mathbb{C})$.

We first show that the condition $[0,p]\cap[0,q]=\{0\}$ does not imply $pq=0$.
Let
\[
p=\begin{pmatrix}1&0\\0&0\end{pmatrix}
\quad\text{and}\quad
q=\frac{1}{2}\begin{pmatrix}1&1\\1&1\end{pmatrix}.
\]
These are the rank-one projections onto $\mathbb{C}(1,0)^{\mathsf{T}}$ and $\mathbb{C}(1,1)^{\mathsf{T}}$, respectively. 
Let $x\in[0,p]\cap[0,q]$. 
We claim that $\operatorname{ran}(x)\subseteq\operatorname{ran}(p)$. Indeed, if $\zeta\in\ker(p)$, then
\[
0\leq\langle x\zeta,\zeta\rangle\leq\langle p\zeta,\zeta\rangle=0.
\]
Since $x\geq0$, it follows that $x^{1/2}\zeta=0$, and hence $x\zeta=0$. 
Thus $\ker(p)\subseteq\ker(x)$, so $\operatorname{ran}(x)\subseteq\operatorname{ran}(p)$. 
The same argument gives $\operatorname{ran}(x)\subseteq\operatorname{ran}(q)$. 
Consequently,
\[
\operatorname{ran}(x)\subseteq \operatorname{ran}(p)\cap\operatorname{ran}(q)=\{0\},
\]
and therefore $x=0$. Hence $[0,p]\cap[0,q]=\{0\}$. 

On the other hand,
\[
pq=\frac{1}{2}\begin{pmatrix}1&1\\0&0\end{pmatrix}\neq0.
\]
In particular, $p\not\perp q$ by Corollary~\ref{cor:ortho}. We next show that $pq=0$ does not imply $p\perp q$.
Let
\[
p=\begin{pmatrix}1&0\\0&0\end{pmatrix} \quad\text{and}\quad
q=\begin{pmatrix}0&0\\0&1\end{pmatrix}.
\]
Clearly, $pq=0$. Proposition~\ref{prop:ab_zero_interval} therefore gives $[0,p]\cap[0,q]=\{0\}$. 
Define
\[
w=\begin{pmatrix}2&\sqrt{3}\\\sqrt{3}&2\end{pmatrix}\in M_2(\mathbb{C})_{\sa}.
\]
Then
\[
w-(p-q)=\begin{pmatrix}1&\sqrt{3}\\\sqrt{3}&3\end{pmatrix}\quad\text{and}\quad
w-(q-p)=\begin{pmatrix}3&\sqrt{3}\\\sqrt{3}&1\end{pmatrix}.
\]
Both matrices have positive diagonal entries and determinant zero, and hence are positive semidefinite. 
Thus $w\geq p-q$ and $w\geq q-p$, or equivalently, $w\in \{p-q,q-p\}^u$. 
However,
\[
w-(p+q)=\begin{pmatrix}1&\sqrt{3}\\\sqrt{3}&1\end{pmatrix}
\]
has determinant $1-3=-2<0$, and is therefore not positive semidefinite. Hence $w\not\geq p+q$, so 
$w\notin \{p+q,-(p+q)\}^u$. It follows that
\[
\{p-q,q-p\}^u\neq \{p+q,-(p+q)\}^u,
\]
and therefore $p\not\perp q$. This also follows from Theorem~\ref{thm:main}, since $pAq\neq\{0\}$. For example, if
\[
e=\begin{pmatrix}0&1\\0&0\end{pmatrix},
\]
then $peq=e\neq0$.
\end{example}

\begin{cor}\label{cor:comm} 
Let $A$ be a \Cstar-algebra. The following conditions are equivalent: 
\begin{enumerate}[(a)] 
\item $A$ is commutative. 
\item $A_{\sa}$ is a vector lattice.
\item $A_{\sa}$ is a weakly pervasive pre-Riesz space.
\item Every pair $u,v\in A_{\sa}$ satisfying $uv=0$ is disjoint. 
\item Every pair of positive contractions $a,b\in A_+$ satisfying $ab=0$ is disjoint. 
\end{enumerate} 
\end{cor}

\begin{bproof} 
$\textup{(a)}\implies\textup{(b)}$: 
It is a classical result that the self-adjoint part $A_{\sa}$ of a $C^{*}$-algebra $A$ is a vector lattice if and only if $A$ is commutative; the forward implication follows immediately from the Gelfand--Naimark representation, while the converse was established by Sherman \cite[Theorems 1 and 2]{Sherman51}.

\par\smallskip\noindent
$\textup{(b)}\implies\textup{(c)}$:
Suppose that $A_{\sa}$ is a vector lattice. For $a,b\in A_+$, we have $[0,a]\cap[0,b]=[0,a\wedge b]$. 
Consequently,
\[
[0,a]\cap[0,b]=\{0\}\quad\Longleftrightarrow\quad a\wedge b=0 \quad\Longleftrightarrow\quad a\perp b.
\]
Thus disjointness and $D$-disjointness coincide on $A_+$, and hence $A_{\sa}$ is weakly pervasive.

\par\smallskip\noindent 
$\textup{(c)}\implies\textup{(d)}$: Suppose that $A_{\sa}$ is weakly pervasive, and let $u,v\in A_{\sa}$ satisfy $uv=0$. Then $u^2v^2=0$, and hence $|u||v|=0$ by continuous functional calculus. By Proposition~\ref{prop:ab_zero_interval}, $[0,|u|]\cap[0,|v|]=\{0\}$. Thus $|u|$ and $|v|$ are $D$-disjoint. Since $A_{\sa}$ is weakly pervasive, we have $|u|\perp|v|$. By Corollary~\ref{cor:absolute_values}, it follows that $u\perp v$.

\par\smallskip\noindent 
$\textup{(d)}\implies\textup{(e)}$: This is immediate since $A_+ \subseteq A_{\sa}$.

\par\smallskip\noindent 
$\textup{(e)}\implies\textup{(d)}$: Let $u,v\in A_{\sa}$ satisfy $uv=0$. If $u=0$ or $v=0$, then $u\perp v$ follows directly from the definition. Suppose that $u\neq0$ and $v\neq0$. Then $u^2v^2=0$, and hence $|u||v|=0$ by continuous functional calculus. Since $\||u|\|=\|u\|>0$ and $\||v|\|=\|v\|>0$, the elements $\widetilde{a}=\|u\|^{-1}|u|$ and $\widetilde{b}=\|v\|^{-1}|v|$ are positive contractions in $A_+$ satisfying $\widetilde{a}\widetilde{b}=0$. Condition~\textup{(e)} gives $\widetilde{a}\perp\widetilde{b}$. By Lemma~\ref{lem:basic}\textup{(b)}, we obtain
\[
|u| = \|u\|\widetilde{a} \perp \|v\|\widetilde{b} = |v|.
\]
Corollary~\ref{cor:absolute_values} then yields $u\perp v$.

\par\smallskip\noindent 
$\textup{(d)}\implies\textup{(a)}$: We use the classical characterization due to Kaplansky that a \Cstar-algebra is commutative if and only if $x^2=0$ implies $x=0$ (see \cite[2.12.21, p.~68]{dix}). Let $x\in A$ satisfy $x^2=0$, and set $a=x^*x$ and $b=xx^*$. Then $a,b\in A_{\sa}$ and 
\[ ab = x^*x^2x^* = 0. \] 
Condition~\textup{(d)} gives $a\perp b$. By Theorem~\ref{thm:main}, we have $aAb=\{0\}$. In particular, $ax^*b=0$. Multiplying on the left by $x$ gives \[ 0 = x(ax^*b) = (xx^*)^3. \] Since $xx^*\geq0$, it follows by the $C^*$ identity that 
\[ 
0 = \|(xx^*)^3\| = \|xx^*\|^3 = \|x\|^6. 
\] 
Thus $x=0$, and so $A$ is commutative. 
\end{bproof}

\begin{remark}
The implication $\textup{(c)}\implies\textup{(a)}$ was proved in \cite[Corollary 5.3]{Jokela26}. 
Condition~\textup{(d)} cannot, in general, be weakened by requiring the implication only for projections. Indeed, consider the noncommutative \Cstar-algebra $A=C_0((0,1],M_2(\mathbb{C}))$ of continuous $M_2(\mathbb{C})$-valued functions on $(0,1]$ vanishing as $t\to 0^+$. For any projection $p\in A$, the function $t\mapsto\|p(t)\|$ is continuous on the connected interval $(0,1]$ and takes values in $\{0,1\}$ (because $\|p\|^2=\|p^*p\|=\|p^2\|=\|p\|$). Since $\|p(t)\|\to 0$ as $t\to 0^+$, we must have $\|p(t)\|=0$ for all $t\in(0,1]$, so $p=0$. Thus $A$ contains no nonzero projections, and the implication $pq=0 \implies p\perp q$ holds vacuously for all projections $p,q\in A$.
\end{remark}

However, if $A$ is a \Cstar-algebra of \emph{real rank zero} (in particular, if $A$ is any von Neumann algebra), then $pq = 0 \implies p \perp q$ on projections does characterize commutativity. %Recall that $A$ has real rank zero if...

\begin{cor}\label{cor:proj} 
Let $A$ be a \Cstar-algebra that is either unital or has real rank zero. The following conditions are equivalent: 
\begin{enumerate}[(a)] 
\item 
Every projection in $A$ belongs to the center $Z(A)$. 
\item 
Every pair of $D$-disjoint projections $p,q\in\operatorname{Proj}(A)$ is disjoint in $A_{\sa}$. 
\item 
Every pair of algebraically orthogonal projections $p,q\in\operatorname{Proj}(A)$ is disjoint in $A_{\sa}$. 
\end{enumerate} 
In particular, if $A$ has real rank zero,  
then these conditions imply that $A$ is commutative. 
\end{cor} 

\begin{bproof} 
$\textup{(a)}\implies\textup{(b)}$: Suppose that every projection in $A$ is central, and let $p,q\in\operatorname{Proj}(A)$ satisfy $[0,p]\cap[0,q]=\{0\}$. Since $p$ and $q$ are central, they commute. Hence $pq$ is a projection. Moreover, \[ 0\leq pq\leq p \quad\text{and}\quad 0\leq pq\leq q. \] Thus $pq\in[0,p]\cap[0,q]$, and consequently $pq=0$. Since $p\in Z(A)$, we have \[ pxq=xpq=0 \] for every $x\in A$. Hence $pAq=\{0\}$, and Theorem~\ref{thm:main} gives $p\perp q$. 

\par\smallskip \noindent 
$\textup{(b)}\implies\textup{(c)}$: Let $p,q\in\operatorname{Proj}(A)$ satisfy $pq=0$. By Proposition~\ref{prop:ab_zero_interval}, we have $[0,p]\cap[0,q]=\{0\}$. Condition~\textup{(b)} therefore gives $p\perp q$.

\par\smallskip \noindent 
$\textup{(c)}\implies\textup{(a)}$: Let $p\in\operatorname{Proj}(A)$. We show that $p\in Z(A)$. Suppose first that $A$ is unital. Then $\mathbf{1}-p\in\operatorname{Proj}(A)$ and $p(\mathbf{1}-p)=0$. Condition~\textup{(c)} gives $p\perp\mathbf{1}-p$, and hence Theorem~\ref{thm:main} gives \[ pA(\mathbf{1}-p)=\{0\}. \] Taking adjoints gives \[ (\mathbf{1}-p)Ap=\{0\}. \] Therefore, for every $x\in A$, \[ px = pxp+px(\mathbf{1}-p) = pxp \] and \[ xp = pxp+(\mathbf{1}-p)xp = pxp. \] Thus $px=xp$ for every $x\in A$, and hence $p\in Z(A)$. 

Assume then that $A$ is nonunital and has real rank zero, and let $p\in\operatorname{Proj}(A)$. All expressions involving $\mathbf{1}$ are understood in the minimal unitization $\widetilde A$. Set \[ B = (\mathbf{1}-p)A(\mathbf{1}-p). \] This is a hereditary \Cstar-subalgebra of $A$ and therefore has an approximate identity $(q_\lambda)$ consisting of projections (c.f. \cite[Theorem~2.6]{BP91}). Since $q_\lambda\in B$, we have $pq_\lambda=0$. Condition~\textup{(c)} and Theorem~\ref{thm:main} therefore give \[ pAq_\lambda=\{0\} \] for every $\lambda$. Fix $x\in A$ and set $y=x(\mathbf{1}-p)$. Then \[ y^*y = (\mathbf{1}-p)x^*x(\mathbf{1}-p) \in B. \] Since $(q_\lambda)$ is an approximate identity for $B$, it follows that $yq_\lambda\to y$ in norm. Indeed, \[ \|y-yq_\lambda\|^2 = \|(\mathbf{1}-q_\lambda)y^*y (\mathbf{1}-q_\lambda)\| \longrightarrow0. \] For every $\lambda$, we have \[ pyq_\lambda = px(\mathbf{1}-p)q_\lambda \in pAq_\lambda = \{0\}. \] Passing to the limit gives \[ px(\mathbf{1}-p)=py=0. \] Applying this conclusion to $x^*$ and taking adjoints yields \[ (\mathbf{1}-p)xp=0. \] Consequently, \[ px=pxp=xp \] for every $x\in A$. Thus $p\in Z(A)$. We have proved the equivalence of \textup{(a)}, \textup{(b)}, and \textup{(c)}.

Finally, suppose that $A$ has real rank zero and satisfies condition~\textup{(b)} or condition~\textup{(c)}. By the equivalence just proved, every projection in $A$ is central. Moreover, by \cite[Theorem~2.6]{BP91}, the real linear span of $\operatorname{Proj}(A)$ is norm-dense in $A_{\sa}$, and hence its complex linear span is norm-dense in $A$. Thus $Z(A)$ contains a norm-dense linear subspace of $A$. Since $Z(A)$ is norm closed, it follows that $Z(A)=A$. Hence, $A$ is commutative. 
%%If $A$ is commutative, then every projection is central, so condition \textup{(a)} holds.
\end{bproof}

\vspace{0.2cm}
\section*{Acknowledgments}
The author acknowledges the use of large language models (Gemini and GPT-5.6) as assistive tools during the preparation of this manuscript. These tools were used for examining proofs, %identifying potentially omitted assumptions or justifications, 
suggesting possible revisions to arguments, checking intermediate calculations, and improving the clarity and LaTeX presentation of the manuscript. All results, proofs, and calculations have been independently verified by the author, who assumes sole responsibility for the content of the paper.


\begin{thebibliography}{9}

\bibitem{vangaans26}
E.~Bilokopytov, M.~Roelands, S.~Tiersma, J.~H.~van~der~Walt, O.~van~Gaans, and M.~Wortel, \emph{Order convergence in ordered vector spaces: compactness, atoms, pervasiveness, and ice cream spaces}, preprint, arXiv:2610.16554 (2026). \url{https://doi.org/10.48550/arXiv.2610.16554}

\bibitem{BP91}
L.~G.~Brown and G.~K.~Pedersen, \emph{$C^*$-algebras of real rank zero}, J. Funct. Anal. \textbf{99} (1991), 131--149.

\bibitem{dix} 
%J.~Dixmier, \emph{$C^*$-algebras}. North-Holland, 1977. %(2.12.21, p. 68)
J.~Dixmier, \emph{$C^*$-Algebras}. North-Holland Mathematical Library, vol.~15, North-Holland, Amsterdam, 1977.

%\bibitem{Jokela25}
%J.~Jokela, \emph{A note on disjointness and discrete elements in partially ordered vector spaces}, preprint, arXiv:2509.11212 (2025).

\bibitem{Jokela26}
J.~Jokela, \emph{A note on disjointness and discrete elements in partially ordered vector spaces}, Indag. Math. (2026). \url{doi:10.1016/j.indag.2026.02.008} %(PII: S0019357726000789).

%\bibitem{kadison}
%R. Kadison, \emph{Order properties of bounded self-adjoint operators}, Proceedings of the American Mathematical Society, (2), (1951), 505--510.

\bibitem{KV}
A.~Kalauch and O.~van~Gaans, \emph{Pre-Riesz Spaces}, De Gruyter, Berlin, 2019.

%\bibitem{Kaplansky1949}
%I.~Kaplansky, \emph{Normed algebras}, Duke Math. J. \textbf{16} (1949), 399--418.

\bibitem{Pedersen2018}
G.~K.~Pedersen, \emph{$C^*$-Algebras and Their Automorphism Groups, 2nd edition}, Academic Press, 2018.
%London Mathematical Society Monographs, Vol.~14, Academic Press, London--New York, 2018.


\bibitem{Sherman51}
S.~Sherman, \emph{Order in operator algebras}, Amer. J. Math. \textbf{73} (1951), 227--232.

%\bibitem{stone}
%M. H. Stone, \emph{A general theory of spectra I}. Proceedings of the National Academy of Sciences, vol. 26 (1940), pp. 280-283. 

%\bibitem{vGK}
%O.~van~Gaans and A.~Kalauch, \emph{Disjointness in partially ordered vector spaces}, Positivity \textbf{10} (2006), 573--589.

%\bibitem{vGK2}
%O.~van~Gaans and A.~Kalauch, \emph{Bands in partially ordered vector spaces}, Positivity \textbf{12} (2008), 553--571.

%\bibitem{Takesaki79}
%M.~Takesaki, \emph{Theory of Operator Algebras I}, Springer-Verlag, New York--Heidelberg, 1979.

\bibitem{vanHaandel1993}
M.~van Haandel, \emph{Completions in Riesz Space Theory}, Ph.D. thesis, Radboud University Nijmegen, 1993.

\end{thebibliography}
\end{document}